\documentclass[12pt]{article}
\usepackage{amssymb}

\input{tcilatex}
\begin{document}

\begin{center}
\textbf{Uniformity and self-neglecting functions:}

\textbf{II. Beurling Regular Variation and the class }$\Gamma $\textbf{\\[0pt%
]
by \\[0pt]
N. H. Bingham and A. J. Ostaszewski}\\[0pt]

\bigskip

\textit{To Paul Embrechts on his 60}$^{\text{th}}$\textit{\ birthday.}

\bigskip
\end{center}

\textbf{Abstract. }Beurling slow variation is generalized to Beurling
regular variation. A Uniform Convergence Theorem, not previously known, is
proved for those functions of this class that are measurable or have the
Baire property. This permits their characterization and representation. This
extends the gamma class of de Haan theory studied earlier.

\textbf{Classification}: 26A03; 33B99, 39B22, 34D05.\newline

\noindent \textbf{1. Introduction}\newline

This paper is a sequel, both to our recent joint paper [BinO11] (Part I\
below), and to the second author's earlier paper [Ost1], in which regular
variation was studied from the viewpoint of topological dynamics in general
and cocycles in particular. It is inspired also by [BinO1]. Our reference
for regular variation is [BinGT] (BGT below).

We begin by setting the context of what we call \textit{Beurling regular
variation} as an extension of the established notion of Beurling slow
variation, and then in Section 2 recall from Part I combinatorial
preliminaries, used there to expand Bloom's analysis of self-neglecting
functions; they are key in enabling us to establish (in Section 4) a
Beurling analogue of the \textit{Uniform Convergence Theorem} (UCT)\ of
Karamata theory. In Section 3 we return to the flow issues raised in this
introduction below; there flow rates, time measures and cocycles are
introduced. Here we discuss the connection between the orbits of the
relevant flows and the Darboux property that plays such a prominent role in
Part I. Incidentally, this explains why the Darboux property is quite
natural in the context of Part I. These ideas prepare the ground for a
Beurling version of the UCT. We deduce a \textit{Characterization Theorem}
in Section 5. Armed with these two theorems, we are able in Section 6 to
establish various \textit{Representation Theorems} for Beurling
regularly-varying functions, but only after a review of Bloom's work on the
representation of self-neglecting functions, from which we glean \textit{%
Smooth Variation Theorems}. We close in Section 7 by commenting on the place
of Karamata theory, and of de Haan's theory of the gamma class, relative to
the new Beurling Theory.

As in Part I the reader should have clearly in mind two isometric
topological groups: the real line under addition with (the Euclidean
topology and) Haar measure Lebesgue measure $dx$, and the positive half-line
under multiplication, with Haar measure $dx/x,$ and metric $d_{W}(x,y)=|\log
y-\log x|.$ (\textquotedblleft W for Weil\textquotedblright , as this
generates the underlying Weil topology of the Haar measure -- for which see
[Halm, \S 62], [Wei]; cf. [BinO7, Th. 6.10], and Part I, \S\ 5.11.) As usual
(see again e.g. Part I, \S\ 5.11) we will move back and forth between these
two as may be convenient, by using their natural isomorphism exp/log. Again
as usual, we work additively in proofs, and multiplicatively in
applications; we use the convention (as in Part I)%
\[
h:=\log f,\qquad k:=\log g.
\]

The new feature Beurling regular variation presents, beyond Karamata regular
variation, is the need to use both addition and multiplication
simultaneously. It is this that makes the affine group $\mathcal{A}ff$ a
natural ingredient here. Recall that on the line the \textit{affine group} $%
x\rightarrow ux+v$ with $u>0$ and $v$ real has (right) Haar measure $u^{-1}du
$ $dv$ (or $(du/u)dv,$ as above) -- see [HewR, IV, (15.29)]. This explains
the presence of the two measure components in the representation of a
Beurling regularly varying function with index $\rho $:%
\begin{equation}
f(x)=d(x)\exp \left( \rho \int_{1}^{x}\frac{u}{\varphi (u)}\frac{du}{u}%
\right) \left( \int_{0}^{x}\frac{e(v)}{\varphi (v)}dv\right) ,  \tag{$\Gamma
_{\rho }$}
\end{equation}%
with $d$ converging to a constant, $e$ smooth and vanishing at infinity, and
integrals initialized at the appropriate group identity ($0$ for $dv$ and $1$
for $du/u)$. Here $u/\varphi (u)$ should be viewed as a \textit{density
function} (for the Haar measure $du/u$ -- unbounded, as $\varphi (x)=o(x)$
in this `Beurling case', but with $\varphi (x)=x$ giving the `Karamata case'
in the limit) -- see Section 3.

Generalizations of Karamata's theory of regular variation (BGT; cf. [Kor]),
rely on a group $G$ acting on a space $X$ in circumstances where one can
interpret `limits to infinity' $x\rightarrow \infty $ in the following
expression:%
\[
g(t):=\lim\nolimits_{x\rightarrow \infty }f(tx)/f(x),\text{ for }t\in G\text{
and }x\in X.
\]%
Here an early treatment is [BajK] followed by [Bal], but a full topological
development dates from the more recent papers [BinO5], [BinO6] and [Ost1] --
see also [BinO10] for an overview. Recall that a group action $A:G\times
X\rightarrow X$ requires two properties:\newline
\noindent (i) \textit{identity}: $A(1_{G},x)=x$ for all $x,$ i.e. $1_{G}=%
\mathrm{id}$, and\newline
\noindent (ii) \textit{associativity}: $A(gh,x)=A(g,A(h,x)),$\newline
with the maps $x\rightarrow g(x):=A(g,x),$ also written $gx,$ often being
homeomorphisms. An action $A$ defines an $A$\textit{-flow} (also referred to
as a $G$-flow), whose orbits are the sets $Ax=\{gx:g\in G\}.$

In fact (i) follows from (ii) for \textit{surjective} $A$ (as $%
A(g,y)=A(1_{G},A(g,y))),$ so we will say that $A\ $is a \textit{pre-action}
if just (i) holds, and then continue to use the notation $g(x):=A(g,x);$ it
is helpful here to think of the corresponding sets $Ax$ as orbits of an $A$%
\textit{-preflow, }using the language of flows and topological dynamics
[Bec].

Beurling's theory of slow variation, introduced in order to generalize the
Wiener Tauberian Theorem ([Kor], Part I), is concerned with consequences of
the equation%
\begin{equation}
f(x+t\varphi (x))/f(x)\rightarrow 1\text{ as }x\rightarrow \infty \text{ }%
\forall t,  \tag{$BSV$}  \label{BSV}
\end{equation}%
for $f:\mathbb{R}\rightarrow \mathbb{R}_{+}$, equivalently,%
\begin{equation}
h(x+t\varphi (x))-h(x)\rightarrow 0\text{ as }x\rightarrow \infty \text{ }%
\forall t,  \tag{$BSV_{+}$}
\end{equation}%
for $h:\mathbb{R}\rightarrow \mathbb{R}$, where $\varphi >0$ itself
satisfies the stronger property:%
\begin{equation}
\varphi (x+t\varphi (x))/\varphi (x)\rightarrow 1\text{ as }x\rightarrow
\infty \text{ locally uniformly in }t\text{.}  \tag{$SN$}
\end{equation}%
Such a function $\varphi $ is said to be \textit{self-neglecting}, briefly $%
\varphi \in SN$. A function $\varphi $ is said to be \textit{Beurling slowly
varying} if $(BSV)$ holds for $f=\varphi $ together with the side condition $%
\varphi (x)=o(x)$. A measurable self-neglecting function is necessarily
Beurling slowly varying in this sense (as $\varphi (x)=o(x);$ see Theorem
4). There are issues surrounding the converse direction, for which see Part
I (as side-conditions are needed for uniformity).

Below we relax the definition of regular variation so that it relies not so
much on group-action but on asymptotic \textquotedblleft cocycle
action\textquotedblright\ associated with a group $G$. This will allow us to
develop a theory of Beurling regular variation analogous to the Karamata
theory, in which the regularly varying functions are those functions $f$
with the Baire property (briefly, \textit{Baire}) or measurable that for
some fixed self-neglecting $\varphi $ possess a \textit{non-zero} limit
function $g$ (not identically zero modulo null/meagre sets) satisfying 
\begin{equation}
f(x+t\varphi (x))/f(x)\rightarrow g(t),\text{ as }x\rightarrow \infty ,\text{
}\forall t,  \tag{$BRV$}
\end{equation}%
(so that $g(0)=1).$ Equivalently,%
\begin{equation}
h(x+t\varphi (x))-h(x)\rightarrow k(t),\text{ as }x\rightarrow \infty ,\text{
}\forall t,  \tag{$BRV_{+}$}
\end{equation}%
(so that $k(0)=0).$ This latter equivalence is non-trivial: it follows from
Theorem 3 below that if $g$ is a non-zero function, then it is in fact
positive. Specializing $(BRV)$ to the sequential format%
\[
g(t)=\lim\nolimits_{n\in \mathbb{N}}f(n+t\varphi (n))/f(n), 
\]%
one sees that the limit function $g$ is Baire/measurable if $f$ is so. We
refer to functions $f$ satisfying $(BRV)$ as (Beurling) $\varphi $\textit{%
-regularly varying}.

This takes us beyond the classical development of such a theory restricted
to the class $\Gamma $ of monotonic functions $f$ satisfying the equation
(BGT\ \S 3.10; de Haan [deH]), and comes on the heels of a recent
breakthrough concerning local uniformity of Beurling slow variation in Part
I. We prove in Theorem 2 below a Uniform Convergence Theorem for
Baire/measurable functions $f$ with non-zero limit $g$, not previously
known, and in Theorem 3 the Characterization Theorem that a Baire/measurable
function $f$ is Beurling $\varphi $-regularly varying with non-zero limit
iff for some $\rho $ one has 
\[
f(x+t\varphi (x))/f(x)\rightarrow e^{\rho t}\text{ }\forall t,
\]%
where $\rho $ is the \textit{Beurling }$\varphi $\textit{-index} of regular
variation. Baire and measurable (positive) functions of this type form the
class $\Gamma _{\rho }(\varphi )$ (cf. Mas [Mas, \S 3.2], Omey [Om], for $f$
measurable; see also [Dom] for the analogous power-wise approach to Karamata
regular variation).

\bigskip

\noindent \textbf{2. Combinatorial preliminaries}\newline

\indent As usual with proofs involving regular variation the nub lies in
infinite combinatorics, to which we now turn. We recall that one can handle
Baire and measurable cases together by working bi-topologically, using the
Euclidean topology in the Baire case (the primary case) and the density
topology in the measure case; see [BinO4], [BinO8], [BinO7]. The \textit{%
negligible }sets are the meagre sets in the Baire case and the null sets in
the measure case; we say that a property holds \textit{quasi everywhere }if
it hold off a negligible set.

We work in the affine group $\mathcal{A}ff$ acting on $(\mathbb{R},+)$ using
the notation 
\[
\gamma _{n}(t)=c_{n}t+z_{n},
\]%
where $c_{n}\rightarrow c_{0}=c>0$ and $z_{n}\rightarrow 0$ as $n\rightarrow
\infty $, as in Theorem 0 below. These are to be viewed as (self-)
homeomorphisms of $\mathbb{R}$ under either the Euclidean topology, or the
Density topology. We recall the following definition from [BinO4] and a
result from [BinO10].

\bigskip

\noindent \textbf{Definition. }A sequence of homeomorphisms $%
h_{n}:X\rightarrow X$ satisfies the \textit{weak category convergence }%
condition (wcc) if:\newline
\indent For any non-meagre open set $U\subseteq X,$ there is a non-meagre
open set $V\subseteq U$ such that for each $k\in N,$%
\[
\bigcap\nolimits_{n\geq k}V\backslash h_{n}^{-1}(V)\text{ is meagre.} 
\]

\noindent \textbf{Theorem 0 (Affine Two-sets Theorem). }\textit{For }$%
c_{n}\rightarrow c>0$ \textit{and }$z_{n}\rightarrow 0,$ \textit{if }$%
cB\subseteq A$ \textit{for }$A,B$ \textit{non-negligible (measurable/Baire),
then for quasi all }$b\in B$ \textit{there exists an infinite set }$\mathbb{M%
}=\mathbb{M}_{b}\subseteq \mathbb{N}$ \textit{such that}%
\[
\{\gamma _{m}(b)=c_{m}b+z_{m}:m\in \mathbb{M}\}\subseteq A.
\]

Below we use only the case $c=1,$ for which $\gamma _{m}$ converge to the
identity (in supremum norm) as a sequence with (wcc), a matter verified in
the Baire case in [BinO6, Th. 6.2] and in the measure case in [BinO3, Cor.
4.1].

\bigskip

\textbf{3. Flows, orbits, cocycles and the Darboux property.}\newline

Our approach is to view Beurling regular variation as a generalization of
Karamata regular variation obtained by replacing the associativity of group
action by a form of asymptotic associativity. To motivate our definition
below, take $X$ and $G$ both to be $(\mathbb{R},\mathbb{+)}$, $\varphi \in SN
$ and consider the map%
\[
T^{\varphi }:(t,x)\mapsto x+t\varphi (x).
\]%
One wants to think of $t$ as representing translation. For fixed $t$ put%
\[
t(x),\text{ or just }tx:=T_{t}^{\varphi }(x)=T^{\varphi }(t,x)=x+t\varphi
(x),
\]%
so that $0(x)=x,$ and so $T^{\varphi }$ is a pre-action. Here we have $%
T_{s+t}^{\varphi }(x)=x+(s+t)\varphi (x),$ so that%
\[
T_{s}^{\varphi }(T_{t}^{\varphi }(x))=x+t\varphi (x)+s\varphi (x+t\varphi
(x))\neq T_{s+t}^{\varphi }(x).
\]%
So $T^{\varphi }:G\times X\rightarrow X$ is not a group action, as
associativity fails. However, just as in a proper flow context, here too one
has a well-defined \textit{flow rate}, or infinitesimal generator, at $x,$
for which see [Bec], [Rud, Ch. 13], [BinO1], cf. [Bal], 
\[
\dot{T}_{0}^{\varphi }x=\left. \frac{d}{dt}T_{t}^{\varphi }x\right\vert
_{t=0}=\lim\nolimits_{t\rightarrow 0}\frac{T_{t}^{\varphi }x-x}{t}=\varphi
(x).
\]%
There is of course an underlying true flow here -- in the measure case,
generated\footnote{%
Positivity is key here; $x=0$ is a fixed-point of the flow $\dot{u}=\varphi
(u)$ when $\varphi (x)=\sqrt{|x|}.$} by $\varphi >0$ (with $1/\varphi $
locally integrable) and described by the system of differential equations
(writing $u_{x}(t)$ for $u(t,x)$) 
\begin{equation}
\dot{u}_{x}(t)=du_{x}(t)/dt=\varphi (u_{x}(t))\text{ with }u_{x}(0)=x,\text{
so that}\quad (t,x)\mapsto u_{x}(t).  \label{ode}
\end{equation}%
(The inverse problem, for $t(u)$ with $t(0)=1,$ has an explicit increasing
integral representation, yielding $u_{x}(t):=u(t+t(x)),$ where $u(t(x))=x,$
as $u$ and $t$ are inverse.) The `differential flow' $\Phi :(t,x)\mapsto
u_{x}(t)$ is continuous in $t$ for each $x.$ As such, $\Phi $ is termed by
Beck a \textit{quasi-flow}.\footnote{%
Beck denotes flows by $\varphi (t,x)$ and uses $f$ where we use $\varphi .$
As we follow the traditional notation of $\varphi $ for self-neglecting
functions, the flow here is denoted $\Phi .$} By contrast `translation
flow', i.e. $(t,x)\mapsto x+t$, being jointly continuous, is a `continuous
flow', briefly a \textit{flow}. It is interesting to note that, by a general
result of Beck (see [Bec] Ch. 4 -- Reparametrization, Th. 4.4.), if the 
\textit{orbits} of $\Phi $ (i.e. the sets $\mathcal{O}(x):=\{\Phi (t,x):t\in 
\mathbb{R}\})$ are continua then, even though $\varphi $ need not be
continuous, there still exists a unique \textit{continuous} `local
time-change' system of maps $t\mapsto f_{x}(t)$ embedding the quasi-flow in
the translation flow, i.e. $\Phi (t,x)=x+f_{x}(t);$ here $f_{x}$ has the 
\textit{cocycle property} (cf. Theorem 1 below),%
\[
f_{x}(s+t)=f_{x}(s)+f_{y}(t)\text{ for }y=x+f_{x}(s),
\]%
and $f_{x}(0)=0$ for all $x.$ This will be the case when $\varphi $ has the
intermediate value property, so here the Darboux property says simply that 
\textit{orbits embed}. Cocycles are thus central to the flow-analysis of
regular variation, central to our earlier index theory of regular variation
[BinO1].

It is this differential flow that the, algebraically much simpler, Beurling
preflow circumvents, working not with the continuous translation function $%
f_{x}(t)$ but $t\varphi (x)$, now only measurable, but with the variables
separated. Nevertheless, the differential equation above is the source of an
immediate interpretation of the integral%
\[
\tau _{x}:=\int_{1}^{x}\frac{du}{\varphi (u)},
\]%
arising in the representation formula ($\Gamma _{\rho }$) for a regularly
varying function $f$, as the metric of \textit{time measure} (in the sense
of Beck -- [Bec, p. 153]). The metric is the \textit{occupation-time measure 
}(cf. BGT\ \S 8.11) of the interval $[1,x]$ under the $\varphi $-generated
flow started at the natural origin of the mutiplicative group $\mathbb{R}_{+}
$. For $\varphi (x)=x,$ the $\varphi $-time measure is Haar measure, and the
associated metric is the Weil (multiplicatively invariant) metric with $%
d_{W}(1,x)=|\log x|,$ as in Section 1. In general, however, the $\varphi $%
-time measure $\mu _{\varphi }$ is obtained from Haar measure via the
density $x/\varphi (x),$ interpretable as a \textit{time-change}
`multiplier' $w(x):=\varphi (x)/x$ (cf. [Bec, 5.41]).

Granted its interpretation, it is only to be expected in ($\Gamma _{\rho }$)
that $\tau _{x}$ multiplies the index $\rho $ describing the asymptotic
behaviour of the function $f.$ The time integral $\tau _{x}$ is in fact
asymptotically equal to the time taken to reach $x$ from the origin under
the Beurling pre-action $T^{\varphi },$ when $\varphi \in SN$, namely $%
x/\varphi (x).$ We hope to return to this matter elsewhere.

Actually, $T^{\varphi }$ is even closer to being an action: it is an \textit{%
asymptotic action} (i.e. asymptotically an action), in view of two
properties critical to the development of regular variation. The first
refers to the dual view of the map $(t,x)\mapsto x(t)=x+t\varphi (x)$ with $x
$ fixed (rather than $t,$ as at the beginning of the section). Here we see
the affine transformation $\alpha _{x}(t)=\varphi (x)t+x.$ This \textit{%
auxiliary group} plays its part through allowing the absorption of a small
\textquotedblleft time\textquotedblright\ variation $t+s$ of $t$ into a
small \textquotedblleft space\textquotedblright\ variation in $x$ involving
a concatenation formula, earlier identified in [BinO1] as a component in the
abstract theory of the index of regular variation.

\bigskip

\noindent \textbf{Lemma 1 (Near-associativity, almost absorption). }\textit{%
For }$\varphi \in SN,$%
\[
T^{\varphi }(t+s,x)=T^{\varphi }(\gamma t,T^{\varphi }(s,x)),\text{ where }%
\gamma =\gamma _{x}^{\varphi }(y_{s})=\varphi (x)/\varphi (y)\text{ and }%
y_{s}:=T^{\varphi }(s,x).
\]%
\textit{and }%
\[
\gamma _{x}^{\varphi }(y_{s})\rightarrow 1\text{ as }s\rightarrow 0.
\]%
\textit{Here} $\gamma $\textit{\ satisfies the concatenation formula}%
\[
\gamma _{x}^{\varphi }(z)=\gamma _{x}^{\varphi }(y)\gamma _{y}^{\varphi
}(z),\qquad \forall x,y,z.
\]%
\textit{Alternatively,}%
\[
T_{t+s}^{\varphi }x=T_{\gamma t}^{\varphi }y_{s},\text{ where }y_{s}=\alpha
_{x}(s)=s\varphi (x)+x,
\]%
\textit{\ equivalently}%
\[
T_{t+s}^{\varphi }x=T_{\beta (t+s)}^{\varphi }\alpha _{x}(s),\text{ for }%
\beta (t)=\gamma _{x}(y_{s})(t-s).
\]%
\textit{\ Proof.} 
\begin{eqnarray*}
T_{t+s}^{\varphi }x &=&x+(t+s)\varphi (x) \\
&=&(x+s\varphi (x))+\frac{\varphi (x)}{\varphi (x+s\varphi (x))}t\varphi
(x+s\varphi (x)).
\end{eqnarray*}%
As for the concatenation formula, one has%
\[
\gamma _{x}^{\varphi }(z)=\frac{\varphi (x)}{\varphi (z)}=\frac{\varphi (x)}{%
\varphi (y)}\cdot \frac{\varphi (y)}{\varphi (z)}=\gamma _{x}^{\varphi
}(y)\gamma _{y}^{\varphi }(z).\qquad \square 
\]

\bigskip

As for the second property of $T^{\varphi },$ recall that for $G$ a group
acting on a second group $X$, a $G$-cocycle on $X$ is a function $\sigma
:G\times X\rightarrow X$ defined by the condition%
\[
\sigma (gh,x)=\sigma (g,hx)\sigma (h,x).
\]%
This definition is already meaningful if a pre-action rather than an action
is defined from $G\times X$ to $X$; so for the purposes of asymptotic
analysis one may capture a weak form of associativity as follows using Lemma
1. For a Banach algebra $X,$ we let $X_{\ast }$ denote its invertible
elements.

\bigskip

\noindent \textbf{Definition.} For $X$ a Banach algebra, given a pre-action $%
T:G\times X\rightarrow X$ (i.e. with $1_{G}x=x$ for all $x,$ where, as
above, $gx:=T(g,x)$), an \textit{asymptotic }$G$\textit{-cocycle} on $X$ is
a map $\sigma :G\times X\rightarrow X_{\ast }$ with the property that for
all $g,h\in G$ and $\varepsilon >0$ there is $r=r(\varepsilon ,g,h)$ such
that for all $x$ with $||x||>r$%
\[
||\sigma (gh,x)-\sigma (g,hx)\sigma (h,x)||_{X}<\varepsilon .
\]%
Say that the cocycle is \textit{locally uniform}, if the inequality holds
uniformly on compact $(g,h)$-sets.

\bigskip

\noindent \textbf{Remark.} Baire and measurable cocycles are studied in
[Ost1] for their uniform boundedness properties. One sees that the Second
and Third Boundedness Theorems proved there hold in the current setting with
asymptotic cocycles replacing cocycles. We now verify that taking $X_{\ast }=%
\mathbb{R}_{+}$, $T=T^{\varphi },$ and the natural cocycle of regular
variation $\sigma ^{f}(t,x):=f(tx)/f(x)$ the above property holds. (More
general contexts are considered in [BinO12] -- Part III.) The case $%
f=\varphi $ comes first; the general case of $\varphi $-regularly varying $f$
must wait till after Theorem 3.

\bigskip

\noindent \textbf{Theorem 1.} \textit{For (positive) }$\varphi \in SN$%
\textit{\ the map }$\sigma ^{\varphi }:\mathbb{R}\times \mathbb{R}%
\rightarrow \mathbb{R}_{+}$%
\[
\sigma ^{\varphi }(t,x):=\varphi (tx)/\varphi (x),\text{ where }%
tx:=T_{t}^{\varphi }(x), 
\]%
\textit{regarded as a map into the Banach algebra }$\mathbb{R}$\textit{\ is
a locally uniform asymptotic }$(\mathbb{R},+)$\textit{-cocycle, i.e. for
every }$\varepsilon >0$\textit{\ and compact set }$K$ \textit{there is }$r$ 
\textit{such that for all }$s,t\in K$\textit{\ and all }$x$ \textit{with } $%
||x||>r$\textit{\ }%
\[
|\sigma ^{\varphi }(s+t,x)-\sigma ^{\varphi }(s,tx)\sigma ^{\varphi
}(t,x)|<\varepsilon , 
\]%
\textit{i.e.}%
\[
\left\vert \frac{\varphi (T_{s+t}^{\varphi }(x))}{\varphi (x)}-\frac{\varphi
(T_{s}^{\varphi }(T_{t}^{\varphi }(x)))}{\varphi (T_{t}^{\varphi }(x))}\cdot 
\frac{\varphi (T_{t}^{\varphi }(x))}{\varphi (x)}\right\vert <\varepsilon , 
\]%
\textit{or}%
\[
|\varphi (T_{s+t}^{\varphi }(x))/\varphi (x)-\varphi (T_{s}^{\varphi
}(T_{t}^{\varphi }(x)))/\varphi (x)|<\varepsilon . 
\]

\bigskip

\noindent \textit{Proof. }Let $\varepsilon >0.$ Given $s,t$ let $I$ be any
open interval with $s+t\in I$.

Pick $\delta >0$ so that the interval $J=(1-\delta ,1+\delta )$ satisfies $%
t+sJ\subseteq I$. Next pick $r$ such that for $x>r$ both%
\[
|\sigma ^{\varphi }(t,x)-1|=|\varphi (x+t\varphi (x))/\varphi (x)-1|<\delta 
\]%
and%
\[
|\sigma ^{\varphi }(u,x)-1|=|\varphi (x+u\varphi (x))/\varphi
(x)-1|<\varepsilon /2,\text{ for all }u\in I.
\]%
In particular%
\[
|\sigma ^{\varphi }(s+t,x)-1|=|\varphi (x+(s+t)\varphi (x))/\varphi
(x)-1|<\varepsilon /2.
\]%
Noting, as in Lemma 1, that 
\[
T_{s}^{\varphi }(T_{t}^{\varphi }x)=(x+t\varphi (x))+s\varphi (x+t\varphi
(x))=x+\varphi (x)\left( t+s\frac{\varphi (x+t\varphi (x))}{\varphi (x)}%
\right) ,
\]%
so that 
\[
w:=t+s\sigma ^{\varphi }(t,x)=t+s\frac{\varphi (x+t\varphi (x))}{\varphi (x)}%
\in t+sJ\subseteq I,
\]%
one has%
\[
|\sigma ^{\varphi }(w,x)-1|<\varepsilon /2,
\]%
i.e.%
\[
|\varphi (T_{s}^{\varphi }(T_{t}^{\varphi }x))/\varphi (x)-1|<\varepsilon /2.
\]%
But%
\begin{eqnarray*}
\sigma ^{\varphi }(w,x) &=&\frac{\varphi (T_{s}^{\varphi }(T_{t}^{\varphi
}x))}{\varphi (x)}=\frac{\varphi (T_{s}^{\varphi }(T_{t}^{\varphi }(x)))}{%
\varphi (T_{t}^{\varphi }(x))}\frac{\varphi (T_{t}^{\varphi }(x))}{\varphi
(x)} \\
&=&\sigma ^{\varphi }(s,T_{t}^{\varphi }(x))\sigma ^{\varphi }(t,x),
\end{eqnarray*}%
so for $x>r$ one has%
\begin{eqnarray*}
&&|\sigma ^{\varphi }(w,x)-\sigma ^{\varphi }(s+t,x)| \\
&\leq &|[\sigma ^{\varphi }(s,T_{t}^{\varphi }(x))\sigma ^{\varphi
}(t,x)]-1|+|\sigma ^{\varphi }(s+t,x)]-1| \\
&\leq &\varepsilon /2+\varepsilon /2=\varepsilon \qquad \square 
\end{eqnarray*}

\noindent \textbf{4. Uniform Convergence Theorem}\newline

\indent We begin with a lemma that yields simplifications later; it implies
a Beurling analogue of the Bounded Equivalence Principle in the Karamata
theory, first noted in [BinO2, Th. 4]. As it shifts attention to the origin,
we call it the Shift Lemma. It has substantially the same statement and
proof as the corresponding Shift Lemma of Part I except that here $h=\log f$
whilst there one has $h=\log \varphi ,$ so that here the difference $%
h(x_{n}+u\varphi (x_{n}))-h(x_{n})$ tends to $k(u)$ rather than to zero. So
we omit the proof. Below \textit{uniform near} a point $u$ means `uniformly
on sequences converging to $u$' and is equivalent to local uniformity at $u$
(i.e. on compact neighbourhoods of $u).$

\bigskip

\noindent \textbf{Lemma 2 (Shift Lemma: uniformity preservation under shift).%
}\textit{\ For any }$u,$\textit{\ convergence in }$(BRV_{+})$\textit{\ is
uniform near }$t=0$\textit{\ iff it is uniform near }$t=u.$

\bigskip

\noindent \textbf{Definition. }Say that $\{u_{n}\}$ with limit $u$ is a 
\textit{witness sequence at }$u$\textit{\ }(for non-uniformity in $h$) if
there are $\varepsilon _{0}>0$ and a divergent sequence $x_{n}$ such that
for $h=\log f$ 
\begin{equation}
|h(x_{n}+u_{n}\varphi (x_{n}))-h(x_{n})|>\varepsilon _{0}\qquad \forall \
n\in \mathbb{N}.  \label{eps-0}
\end{equation}

\bigskip

\noindent \textbf{Theorem 2 (UCT for }$\varphi $\textbf{-regular variation).}
\textit{For }$\varphi \in SN$\textit{, if }$f$ \textit{has the Baire
property (or is measurable)} \textit{and satisfies }($BRV$) \textit{with
limit }$g$ \textit{strictly positive on a non-negligible set, then }$f$%
\textit{\ is locally uniformly }$\varphi $-$RV$\textit{.}

\bigskip

\noindent \textit{Proof.} Suppose otherwise. We modify a related proof from
Part I (concerned there with the special case of $\varphi $ itself) in two
significant details. In the first place, we will need to work relative to
the set $S:=\{s>0:g(s)>0\},$ (\textquotedblleft S for
support\textquotedblright ), so that $k(s)=\log g(s)$ is well-defined on $S.$
Now $S$ is Baire/measurable; as $S$ is assumed non-negligible, by passing to
a Baire/measurable subset of $S$ if necessary, we may assume w.l.o.g. that
the restriction $k|S$ is continuous on $S,$ by [Kur, \S 28] in the Baire
case (cf. [BinO6] Th. 11.8) and Luzin's Theorem in the measure case ([Oxt],
Ch. 8, cf. [BinO7]).

Let $u_{n}$ be a witness sequence for the non-uniformity of $h$ so, for some 
$x_{n}\rightarrow \infty $ and $\varepsilon _{0}>0$ one has (\ref{eps-0}).
By the Shift Lemma (Lemma 2) we may assume that $u=0.$ So we will write $%
z_{n}$ for $u_{n}.$ As $\varphi $ is self-neglecting we have%
\begin{equation}
c_{n}:={\varphi (x_{n}+z}_{n}{\varphi (x_{n}))}/{\varphi
(x_{n})\longrightarrow 1.}  \label{lim}
\end{equation}%
Write $\gamma _{n}(s):=c_{n}s+z_{n}$ and $y_{n}:=x_{n}+z_{n}\varphi (x_{n}).$
Then $y_{n}=x_{n}(1+z_{n}\varphi (x_{n})/x_{n})\rightarrow \infty ,$ and as $%
k(0)=0,$%
\begin{equation}
\left\vert h{(y}_{n}{)-h(x_{n})}\right\vert \geq \varepsilon _{0}.
\label{eps-0'}
\end{equation}

\indent Now take $\eta =\varepsilon _{0}/4$ and for $x=\{x_{n}\}$, working
in $S,$ put 
\[
V_{n}^{x}(\eta ):=\{s\in S:|h(x_{n}+s\varphi (x_{n}))-h(x_{n})-k\left(
s\right) |\leq \eta \},\text{ }H_{k}^{x}(\eta ):=\bigcap\nolimits_{n\geq
k}V_{n}^{x}(\eta ),
\]%
and likewise for $y=\{y_{n}\}.$ These are Baire sets, and 
\begin{equation}
S=\bigcup\nolimits_{k}H_{k}^{x}(\eta )=\bigcup\nolimits_{k}H_{k}^{y}(\eta ),
\label{cov}
\end{equation}%
as $h\in BRV_{+}$. The increasing sequence of sets $\{H_{k}^{x}(\eta )\}$
covers $S.$ So for some $k$ the set $H_{k}^{x}(\eta )$ is non-negligible. As 
$H_{k}^{x}(\eta )$ is non-negligible, by (\ref{cov}), for some $l$ the set%
\[
B:=H_{k}^{x}(\eta ))\cap H_{l}^{y}(\eta )
\]%
is also non-negligible. Taking $A:=H_{k}^{x}(\eta ),$ one has $B\subseteq
H_{l}^{y}(\eta )$ and $B\subseteq A$ with $A,B\ $non-negligible. Applying
Theorem 0 to the maps $\gamma _{n}(s)=c_{n}s+z_{n}$ with $c=\lim_{n}c_{n}=1,$
there exist $b\in B$ and an infinite set\textit{\ }$\mathbb{M}$ such that%
\[
\{c_{m}b+z_{m}:m\in \mathbb{M}\}\subseteq A=H_{k}^{x}(\eta ).
\]%
That is, as $B\subseteq H_{l}^{y}(\eta ),$ there exist $t\in H_{l}^{y}(\eta )
$ and an infinite $\mathbb{M}_{t}$ such that%
\[
\{\gamma _{m}(t)=c_{m}t+z_{m}:m\in \mathbb{M}_{t}\}\subseteq H_{k}^{x}(\eta
).
\]%
In particular, for this $t$ and $m\in \mathbb{M}_{t}$ with $m>k,l$ one has 
\[
t\in V_{m}^{y}(\eta )\hbox{ and }\gamma _{m}(t)\in V_{m}^{x}(\eta ).
\]%
As $t\in S$ and $\gamma _{m}(t)\in S$ (a second, critical, detail), we have
by continuity of $k|S$ at $t,$ since $\gamma _{m}(t)\rightarrow t,$ that for
all $m$ large enough%
\begin{equation}
|k(t)-k(\gamma _{m}(t))|\leq \eta .  \label{cont}
\end{equation}%
Fix such an $m.$ As $\gamma _{m}(t)\in V_{m}^{x}(\eta ),$%
\begin{equation}
\left\vert h{(x_{m}+\gamma }_{m}({t)\varphi (x_{m}))}-h{(x_{m})-k(\gamma }%
_{m}(t))\right\vert \leq \eta .  \label{*}
\end{equation}%
But $\gamma _{m}(t)=c_{m}t+z_{m}=z_{m}+t\varphi (y_{m})/\varphi (x_{m}),$ so 
\[
x_{m}+\gamma _{m}(t)\varphi (x_{m})=x_{m}+z_{m}\varphi (x_{m})+t\varphi
(y_{m})=y_{m}+t\varphi (y_{m}),
\]%
`absorbing' the affine shift ${\gamma }_{m}(${$t$}${)}$ into $y.$ So, by (%
\ref{*}), 
\[
\left\vert h{(y_{m}+t\varphi (y_{m}))}-h{(x_{m})-k(\gamma }%
_{m}(t))\right\vert \leq \eta .
\]%
But $t\in V_{m}^{y}(\eta ),$ so%
\[
\left\vert h{(y_{m}+t\varphi (y_{m}))}-h{(y_{m})-k(t))}\right\vert \leq \eta
.
\]%
By the triangle inequality, combining the last two inequalities with (\ref%
{cont}), 
\begin{eqnarray*}
&&|h(y_{m})-h(x_{m})| \\
&\leq &|h(y_{m}+t\varphi (y_{m}))-h(y_{m})-k(t)|+|k(t)-k(\gamma
_{m}(t))|+|h(y_{m}+t\varphi (y_{m}))-h(x_{m})-k(\gamma _{m}(t))| \\
&\leq &3\eta <\varepsilon _{0},
\end{eqnarray*}%
contradicting (\ref{eps-0'}). \U{25a1}

\bigskip

\noindent \textbf{5. Characterization Theorem}\newline

\indent We may now deduce the characterization theorem which implies in
particular that the support set $S$ of the last proof is in fact all of $%
\mathbb{R}.$

\bigskip

\noindent \textbf{Theorem 3 (Characterization Theorem).} \textit{For }$%
\varphi \in SN,$\textit{\ if }$f>0$\textit{\ is }$\varphi $\textit{%
-regularly varying and Baire/measurable and satisfies}%
\[
f(x+t\varphi (x))/f(x)\rightarrow g(t),\text{  }\forall t,
\]%
\textit{with non-zero limit, i.e. }$g>0$\textit{\ on a non-negligible set,
then for some }$\rho $\textit{\ (the index of }$\varphi $\textit{-regular
variation) one has}%
\[
g(t)=e^{\rho t}.
\]

\bigskip

\noindent \textit{Proof.} Proceed as in Theorem 1: writing $y:=x+s\varphi
(x),$ and recalling from Th. 1 the notation $\gamma =\varphi (x)/\varphi
(y), $ one has%
\begin{equation}
h(x+(s+t)\varphi (x))-h(x)=[h(y+t\gamma \varphi (y))-h(y)]+[h(y)-h(x)].
\label{id+}
\end{equation}

Fix $s$ and $t\in \mathbb{R}$; passing to limits and using uniformity (by
Theorem 2), we have%
\begin{equation}
k(s+t)=k(t)+k(s),  \tag{$CFE$}
\end{equation}%
since $\gamma =\varphi (x)/\varphi (y)\rightarrow 1.$ This is the Cauchy
functional equation; as is well-known, for $k$ Baire/measurable (see Banach
[Ban, Ch. I, \S 3, Th. 4] and Mehdi [Meh] for the Baire case, [Kucz, 9.4.2]
for the measure case, and [BinO8] for an up-to-date discussion) this implies 
$k(x)=\rho x$ for some $\rho \in \mathbb{R}$, and so $g(x)=e^{\rho x}.$ $%
\square $

\bigskip

\noindent \textbf{Remark. }The conclusion that $k(x)=\rho x$ ($\forall x)$
for some $\rho $ tells us that in fact $g>0$ everywhere, which in turn
implies the cocycle property below. (If we assumed that $g>0$ everywhere, we
could argue more naturally, and more nearly as in the Karamata theory, by
establishing the cocycle property first and from it deducing the
Characterization Theorem.)

\bigskip

As an immediate corollary, we now have an extension to Theorem 1:

\bigskip

\noindent \textbf{Corollary 1 (Cocycle property).} \textit{For }$\varphi \in
SN,$\textit{\ if }$f>0$\textit{\ is }$\varphi $\textit{-regularly varying
and Baire/measurable and satisfies}%
\[
f(x+t\varphi (x))/f(x)\rightarrow g(t),\text{  }\forall t,
\]%
\textit{with non-zero limit, i.e. }$g>0$\textit{\ on a non-negligible set,
then}%
\[
\sigma ^{f}(t,x):=f(x+t\varphi (x))/f(x)
\]%
\textit{is a locally uniform asymptotic cocycle.}

\bigskip

\noindent \textit{Proof.} With the notation of Theorem 3, rewrite (\ref{id+}%
) as%
\begin{equation}
\frac{f(x+(s+t)\varphi (x))}{f(x)}=\frac{f(y+t\gamma \varphi (y))}{f(y)}%
\cdot \frac{f(x+s\varphi (x))}{f(x)},  \label{id}
\end{equation}%
where by Theorem 3 both ratios on the right-hand side have non-zero limits $%
g(t)$ and $g(s),$ as $x$ (and so $y)$ tend to infinity. Given $\varepsilon >0
$ it now follows from (\ref{id}) using ($CFE$), and uniformity in compact
neighbourhoods of $s,$ $t$ and $s+t$ (by Theorem 2), that for all large
enough $x$ 
\[
\left\vert \frac{f(x+(s+t)\varphi (x))}{f(x)}-\frac{f(x+t\varphi (x))}{f(x)}%
\cdot \frac{f(x+s\varphi (x))}{f(x)}\right\vert <\varepsilon ,
\]%
so that $\sigma ^{f}$ is a locally uniform asymptotic cocycle. $\square $

\bigskip

\noindent \textbf{6. Smooth Variation and Representation Theorems}\newline

Before we derive a representation theorem for Baire/measurable Beurling
regularly varying functions, we need to link the Baire case to the measure
case. Recall the \textit{Beck iteration} of $\gamma (x):=T_{1}^{\varphi
}(x)=x+\varphi (x)$ (so that $\gamma _{n+1}(x)=\gamma (\gamma _{n}(x))$ with 
$\gamma _{1}=\gamma ,$ for which see [Bec, 1.64] in the context of bounding
a flow) and Bloom's result for\textit{\ }$\varphi \in SN$ concerning the
sequence $x_{n+1}=\gamma _{n}(x_{1}),$ i.e. $x_{n+1}:=x_{n}+\varphi (x_{n}),$
that for all $x_{1}$ large enough one has $x_{n}\rightarrow \infty ,$ i.e.
the sequence gives a \textit{Bloom partition} of $\mathbb{R}_{+}$ (see
[Blo], or BGT \S 2.11). We next need to recall a construction due to Bloom
in detail as we need a slight amendment.

\bigskip

\noindent \textbf{Lemma 4 (Interpolation Lemma).} \textit{For }$\varphi \in
SN,$\textit{\ set }$x_{n+1}:=x_{n}+\varphi (x_{n}),$\textit{\ with }$x_{1}$%
\textit{\ large enough so that }$x_{n}\rightarrow \infty .$ \textit{Put }$%
x_{0}=0.$

\textit{For }$\psi >0$ \textit{a }$\varphi $\textit{-slowly varying
function, there exists a continuously differentiable function }$\phi >0$%
\textit{\ such that}\newline
i) $\phi (x_{n})=\psi (x_{n})$ for $n=0,1,2,...,$\newline
ii) $\phi (x)$ \textit{lies between }$\psi (x_{n})$\textit{\ and }$\psi
(x_{n-1})$\textit{\ for }$x$\textit{\ between }$x_{n-1}$\textit{\ and }$x_{n}
$\textit{\ for }$n=1,2,...,$\newline
iii) $|\phi ^{\prime }(x)|\leq 2|\psi (x_{n})-\psi (x_{n-1})|/\varphi
(x_{n-1}),$ \textit{for }$x$\textit{\ between }$x_{n-1}$\textit{\ and }$x_{n}
$\textit{\ for }$n=1,2,....$

\bigskip

\bigskip

\noindent \textbf{Proof.} Proceed as in [Blo] or BGT \S 2.11; we omit the
details. $\square $

\bigskip

\noindent \textbf{Definition. }Call any function $\phi $ with the properties
(i)-(iii) a ($\varphi $-) \textit{interpolating function} for $\psi .$

\bigskip

We now deduce an extension of the Bloom-Shea Representation Theorem in the
form of a Smooth Beurling Variation Theorem (for smooth variation, see BGT 
\S 2.1.9, following Balkema et al. [BalGdH]). Indeed, the special case $\psi
=\varphi $ is included here. Our proof is a variant on Bloom's. It will be
convenient to introduce:

\bigskip

\noindent \textbf{Definition (Asymptotic equivalence). }For $\varphi ,\phi
>0 $ write $\varphi \sim \phi $ if $\varphi (x)/\phi (x)\rightarrow 1$ as $%
x\rightarrow \infty $. If $\phi \in \mathcal{C}^{1},$ say that $\phi $ is a 
\textit{smooth representation} of $\varphi .$

\bigskip

\noindent \textbf{Theorem 4 (Smooth Beurling Variation).} \textit{For }$%
\varphi \in SN$ \textit{and }$\psi $ \textit{a }$\varphi $\textit{-slowly
varying function, if }$\phi $\textit{\ is any continuously differentiable
function interpolating }$\psi $\textit{\ w.r.t. }$\varphi $ \textit{as in
Lemma 4, then}%
\[
\psi (x)=c(x)\phi (x),\text{ for some positive }c(.)\rightarrow 1,\text{
i.e. }\psi \sim \phi \in \mathcal{C}^{1},
\]%
\textit{\ so that }$\phi $\textit{\ is }$\varphi $\textit{-slowly varying,
and also}%
\[
\left\vert \varphi (x)\frac{\phi ^{\prime }(x)}{\phi (x)}\right\vert
\rightarrow 0.
\]%
\textit{This yields the representation}%
\[
\psi (x)=c(x)\phi (x)=c(x)\exp \left( \int_{0}^{x}\frac{e(u)}{\varphi (u)}%
du\right) ,\text{ for }e\in \mathcal{C}^{1}\text{ with }e\rightarrow 0.
\]%
\textit{Moreover, if }$\psi $ \textit{is Baire/measurable, then so is }$c(x)$%
.

\textit{Furthermore, if }$\psi \in SN,$\textit{\ in particular when }$\psi
=\varphi ,$ \textit{then:}

(i)\textit{\ }$\phi (x)$\textit{\ is self-neglecting and }$\psi \sim \phi
\sim \int_{0}^{x}e(u)du$ \textit{for some continuous }$e$\textit{\ with }$%
e\rightarrow 0;$

(ii) \textit{both }$\phi (x)/x\rightarrow 0$ \textit{and }$\psi
(x)/x\rightarrow 0,$ \textit{as }$x\rightarrow \infty ;$

(iii) \textit{if }$f$ \textit{is }$\psi $\textit{-regularly varying with
index }$\rho $\textit{, then }$f$\textit{\ is }$\phi $\textit{-regularly
varying with index }$\rho $ \textit{with }$\psi \sim \phi \in \mathcal{C}%
^{1}\cap SN.$

\bigskip

\noindent \textit{Proof.} Note that for any $y_{n}$ between $x_{n}$ and $%
x_{n+1}$ one has $\varphi (y_{n})/\varphi (x_{n})\rightarrow 1,$ and so also 
$\psi (y_{n})/\psi (x_{n})=[\psi (y_{n})/\varphi (x_{n})][\varphi
(x_{n})/\varphi (y_{n})][\varphi (y_{n})/\psi (x_{n})]\rightarrow 1;$ indeed 
$y_{n}=x_{n}+t_{n}\varphi (x_{n})$ for some $t_{n}\in \lbrack 0,1],$ and so
the result follows from local uniformity in $\psi $ and $\varphi $ and
because $\psi $ is $\varphi $-slowly varying. This implies first that, if
(say) $\psi (x_{n})\leq \psi (x_{n+1})$, then for $\phi $ as in the
statement of the theorem, 
\[
\frac{\psi (x_{n})}{\psi (x_{n+1})}\leq \frac{\phi (y_{n})}{\psi (y_{n})}%
\leq \frac{\psi (x_{n+1})}{\psi (y_{n+1})},
\]%
and so $\phi (y_{n})/\psi (y_{n})\rightarrow 1;$ similarly for $\psi
(x_{n+1})\leq \psi (x_{n}).$ So $\phi (x)/\psi (x)\rightarrow 1,$ as $%
x\rightarrow \infty .$ So, by Lemma 4 (i) and as $\psi $ is $\varphi $%
-slowly varying,%
\[
\frac{\phi (y_{n})}{\varphi (x_{n})}=\frac{\phi (y_{n})}{\psi (y_{n})}\cdot 
\frac{\psi (y_{n})}{\varphi (x_{n})}=\frac{\psi (y_{n})}{\psi (y_{n})}\cdot 
\frac{\psi (y_{n})}{\varphi (x_{n})}\rightarrow 1,
\]%
i.e. $\phi $ is $\varphi $-slowly varying. Furthermore, by Lemma 4 (iii) and
since $\phi $ and $\psi $ are $\varphi $-slowly varying,%
\[
\left\vert \frac{\varphi (x)}{\phi (x)}\phi ^{\prime }(x)\right\vert \leq
2\left( \frac{\psi (x_{n+1})}{\psi (x_{n})}-1\right) (\psi (x_{n})/\varphi
(x_{n}))\cdot (\varphi (x)/\phi (x))\rightarrow 0.
\]%
Take $c(x):=\psi (x)/\phi (x);$ then $\lim_{x\rightarrow \infty }c(x)=1,$
and re-arranging one has%
\[
\psi (x)=c(x)\phi (x),
\]%
as asserted.

From here we have, setting $e(u):=\varphi (u)\phi ^{\prime }(u)/\phi
(u)\rightarrow 0$ and noting that $e(x)/\varphi (x)\geq 0$ is the derivative
of $\log \phi (x),$%
\[
\psi (x)=c(x)\phi (x)=c(x)\exp \int_{0}^{x}\frac{e(u)}{\varphi (u)}du. 
\]%
Conversely, such a representation yields slow $\varphi $-variation: by the
Mean Value Theorem, for any $t$ there is $s=s(x)\in \lbrack 0,t]$ such that%
\[
\int_{x}^{x+t\varphi (x)}\frac{e(u)}{\varphi (u)}du=\frac{e(x+s\varphi (x))}{%
\varphi (x+s\varphi (x))}\varphi (x+s\varphi (x)), 
\]%
which tends to $0$ uniformly in $t$ as $x\rightarrow \infty ,$ since $%
e(.)\rightarrow 0$ and $\varphi \in SN.$

Now suppose additionally that $\psi \in SN$ (e.g. if $\psi =\varphi $). We
check that then $\phi \in SN$. Indeed, suppose that $u_{n}\rightarrow u;$
then as $\phi (x_{n})=\psi (x_{n}),$ and since $\psi \in SN$, writing $%
y_{n}=x_{n}+u_{n}\phi (x_{n})$ one has 
\[
\frac{\phi (y_{n})}{\phi (x_{n})}=\frac{\psi (x_{n}+u_{n}\phi
(x_{n}))/c(y_{n})}{\psi (x_{n})/c(x_{n})}=\frac{\psi (x_{n}+u_{n}\psi
(x_{n}))/c(y_{n})}{\psi (x_{n})/c(x_{n})}\rightarrow 1, 
\]%
as asserted in (i). Given this $\phi $ apply Lemma 4 with $\phi $ for $%
\varphi $ and $\psi =\varphi =\phi $ to yield a further smooth representing
function $\bar{\phi}\sim \phi $. Then by Lemma 4 (iii) we have for a
corresponding sequence $\bar{x}_{n}$%
\[
|\bar{\phi}^{\prime }(x)|\leq 2|\phi (\bar{x}_{n})-\phi (\bar{x}%
_{n-1})|/\phi (\bar{x}_{n-1})=2|\phi (\bar{x}_{n})/\phi (\bar{x}%
_{n-1})-1|\rightarrow 0. 
\]%
So taking $e(x)=\bar{\phi}^{\prime }(x)$ one has $\lim_{x\rightarrow \infty
}e(x)=0$ and for $\bar{c}(x):=\phi (x)/\bar{\phi}(x),$ one has again $%
\lim_{x\rightarrow \infty }\bar{c}(x)=1.$ So integrating $\bar{\phi}^{\prime
},$ one has%
\[
\phi (x)=\bar{c}(x)\bar{\phi}=\bar{c}(x)\int_{0}^{x}e(u)du\text{ with }%
e(u)\rightarrow 0. 
\]%
From this integral representation, one can check that $\phi $ is
self-neglecting (as in [Blo], BGT \S 2.11).

As to (ii), we first prove this for $\varphi $ itself. So specializing (i)
to $\psi =\varphi ,$ write 
\[
\varphi (x)\sim \int_{0}^{x}e(u)du\text{ with }e(u)\rightarrow 0, 
\]%
we deduce immediately that $\varphi (x)/x\sim
\int_{0}^{x}e(u)du/x\rightarrow 0$.

We now use the fact that $\varphi (x)/x\rightarrow 0$ to consider a general $%
\psi \in SN$ that is $\varphi $-slowly varying. Take $\psi \sim \phi \in 
\mathcal{C}^{1}.$Take $a_{n}=\psi (x_{n})/\varphi (x_{n}),$ $b_{n}=\varphi
(x_{n})/\varphi (x_{n-1})>0,$ so that $a_{n}\rightarrow 1$ and $%
b_{n}\rightarrow 1.$ Put $z_{n}:=\varphi (x_{n})/x_{n}>0,$ so that $%
z_{n}\rightarrow 0,$ as just shown. Now one has by Lemma 4 (i) that 
\begin{eqnarray*}
\frac{\phi (x_{n})}{x_{n}} &=&\frac{\psi (x_{n})}{x_{n-1}+\varphi (x_{n})}=%
\frac{a_{n}}{1+x_{n-1}/\varphi (x_{n})} \\
&=&\frac{a_{n}}{1+(1/z_{n-1})\varphi (x_{n-1})/\varphi (x_{n})}=\frac{a_{n}}{%
1+1/(z_{n-1}b_{n})}\rightarrow 0,
\end{eqnarray*}%
as required.

\bigskip

As to (iii) for $\psi \in SN,$ if $\psi \sim \phi \in \mathcal{C}^{1}$, then 
$\phi \in SN$ (by (ii)). So as $\psi (x)/\phi (x)\rightarrow 1,$ by Theorem 2%
\[
\lim\nolimits_{x\rightarrow \infty }\frac{f(x+t\psi (x))}{f(x)}%
=\lim\nolimits_{x\rightarrow \infty }\frac{f(x+t[\psi (x)/\phi (x)]\phi (x))%
}{f(x)}. 
\]%
That is, $f$ is $\phi $-regularly varying. $\square $

\bigskip

We have just seen that self-neglecting functions are necessarily $o(x).$ We
now see that, for $\varphi \in SN$, a $\varphi $-slowly-varying function is
also $SN$ if it is $o(x).$

\bigskip

\noindent \textbf{Theorem 5.} \textit{For }$\varphi \in SN,$\textit{\ if }$%
\psi >0$\textit{\ is }$\varphi $\textit{-slowly varying and }$\psi (x)=o(x)$%
\textit{, then }$\psi $\textit{\ is SN, and so has a representation}%
\[
\psi \sim \int_{0}^{x}e(u)du\text{ with continuous }e(.)\rightarrow 0. 
\]

\bigskip

\noindent \textit{Proof.}\textbf{\ }Since self-neglect is preserved under
asymptotic equivalence, without loss of generality we may assume that $\psi $
is smooth. Now $\psi (x)/\varphi (x)\rightarrow 1$ (by definition), so for
fixed $u,$ $u[\psi (x)/\varphi (x)]\rightarrow u.$ For $\psi $ a $\varphi $%
-slowly varying function, by the UCT for $\varphi $-regular variation%
\[
\psi (x+t\varphi (x))/\varphi (x)\rightarrow 1,\text{ loc. unif. in }u,
\]%
as $\psi $ is measurable. So in particular, 
\[
\frac{\psi (x+t\psi (x))}{\psi (x)}=\frac{\psi (x+t[\psi (x)/\varphi
(x)]\varphi (x))}{\varphi (x)}\frac{\varphi (x)}{\psi (x)}\rightarrow 1.
\]%
That is, $\psi $ is BSV, since $\psi (x)=o(x).$ But $\psi $ is continuous,
so by Bloom's theorem ([Blo]) $\psi \in SN.$ $\square $

\bigskip

\noindent \textbf{Lemma 5. }\textit{For measurable} $\varphi \in SN,\ $%
\textit{the function}%
\[
f_{\rho }(x):=\exp \left( \rho \int_{1}^{x}\frac{du}{\varphi (u)}\right) 
\]%
\textit{\ is }$\varphi $\textit{-regularly varying with index }$\rho .$

\bigskip

\noindent \textit{Proof. }With $h_{\rho }=\log f_{\rho }$\textit{\ }one has%
\textit{\ }that\textit{\ }%
\begin{eqnarray*}
&&h_{\rho }(x+t\varphi (x))-h_{\rho }(x)-\rho t \\
&=&\rho \int_{x}^{x+t\varphi (x)}\frac{du}{\varphi (u)}-\rho t=\rho
\int_{x}^{x+t\varphi (x)}\left( \frac{\varphi (x)}{\varphi (u)}-1\right) 
\frac{du}{\varphi (x)} \\
&=&\rho \int_{0}^{t}\left( \frac{\varphi (x)}{\varphi (x+v\varphi (x))}%
-1\right) dv=o(1).\qquad \square 
\end{eqnarray*}

\bigskip

We may now establish our main result with $f_{\rho }$ as above.

\bigskip

\noindent \textbf{Theorem 6 (Beurling Representation Theorem).} \textit{For }%
$\varphi \in SN$\textit{\ with }$\varphi $ \textit{Baire/measurable
eventually bounded away from 0, and }$f$\textit{\ measurable and }$\varphi $%
\textit{-regularly varying: for some }$\rho \in \mathbb{R}$\textit{\ and }$%
\varphi $\textit{-slowly varying function }$\tilde{f},$\textit{\ one has}%
\[
f(x)=f_{\rho }(x)\tilde{f}(x)=\exp \left( \rho \int_{1}^{x}\frac{du}{\varphi
(u)}\right) \tilde{f}(x).
\]%
\textit{Any function of this form is }$\varphi $\textit{-regularly varying
with index }$\rho .$

\textit{So }$f\sim f_{\rho }\phi $ \textit{for some smooth representation }$%
\phi $ \textit{of }$\tilde{f}.$

\bigskip

\noindent \textit{Proof.} By Theorem 4(iii), we may assume that $\varphi $
is smooth. Choose $\rho $ as in Theorem 3 and, referring to the flow rate $%
\varphi (x)>0$ at $x,$ put%
\[
\tilde{h}(x):=h(x)-\rho \int_{1}^{x}\frac{du}{\varphi (u)}, 
\]%
where $h=\log f.$ So $\tilde{h}(x)$ is Baire/measurable as $h$ is.

By Theorem 2 (UCT), locally uniformly in $t$ one has a `reduction' formula
for $\tilde{h}:$ 
\[
\tilde{h}(x+t\varphi (x))-\tilde{h}(x)-\rho t+\rho \int_{x}^{x+t\varphi (x)}%
\frac{du}{\varphi (u)}=h(x+t\varphi (x))-h(x)-\rho t=o(1).
\]%
So substituting $u=x+v\varphi (x)$ in the last step,%
\begin{eqnarray*}
\tilde{h}(x+t\varphi (x))-\tilde{h}(x) &=&\rho t-\rho \int_{x}^{x+t\varphi
(x)}\frac{du}{\varphi (u)}+o(1) \\
&=&\rho \int_{x}^{x+t\varphi (x)}\left( 1-\frac{\varphi (x)}{\varphi (u)}%
\right) \frac{du}{\varphi (x)}+o(1) \\
&=&\rho \int_{0}^{t}\left( 1-\frac{\varphi (x)}{\varphi (x+v\varphi (x))}%
\right) dv+o(1) \\
&=&o(1),
\end{eqnarray*}%
and the convergence under the integral here is locally uniform in $t$ since $%
\varphi \in SN.$ So $\exp (\tilde{h})$ is Beurling $\varphi $-slowly
varying. The converse was established in Lemma 5. The remaining assertion
follows from Theorem 4. $\square $

\bigskip

As a second corollary of Theorems 2 and 3 and of the de Bruijn-Karamata
Representation Theorem (see BGT, Ths. 1.3.1 and 1.3.3 and the recent
generalization [BinO10]), we deduce a Representation Theorem for Beurling
regular variation which extends previous results concerned with the class $%
\Gamma $ -- see BGT, Th. 3.10.6. We need the following result, which is
similar to Bloom's Th. 4 except that we use regularity of $\varphi $ rather
than assume conditions on convergence rates.

\bigskip

\noindent \textbf{Lemma 6 (Karamata slow variation).} \textit{If }$\varphi
\in SN$\textit{\ with }$\varphi $ \textit{Baire/measurable eventually
bounded away from 0, then }%
\[
\varphi (x+v)\left/ \varphi (x)\right. \rightarrow 1\text{ as }x\rightarrow
\infty ,\text{ locally uniformly in }v.
\]

\bigskip

\noindent \textit{Proof.} W.l.o.g suppose that $0<K<\varphi (x)$ for all $x$%
. Fix $v;$ then $0\leq |v|/\varphi (x)\leq |v|K^{-1}$ for all $x.$ Let $%
\varepsilon >0.$ Since $\varphi \in SN,$ there is $X=X(\varepsilon ,v)$ such
that%
\begin{equation}
\left\vert \varphi (x+t\varphi (x))\left/ \varphi (x)\right. -1\right\vert
<\varepsilon ,  \label{eps}
\end{equation}%
for all $|t|\leq |v|K^{-1}$ and all $x\geq X.$ So in particular, for $x\geq X
$ and $t:=v/\varphi (x),$ since $|t|\leq vK^{-1},$ substitution in (\ref{eps}%
) yields%
\[
\left\vert \varphi (x+v)\left/ \varphi (x)\right. -1\right\vert <\varepsilon
,
\]%
for $x\geq X.$ This shows that for each $v\in \mathbb{R}$ 
\[
\varphi (x+v)\left/ \varphi (x)\right. \rightarrow 1.
\]%
So $\log \varphi $ is Karamata slowly varying in the additive sense; being
Baire/measurable, by the UCT\ of additive Karamata theory, convergence to
the limit for $\log \varphi ,$ and so convergence for $\varphi $ as above,
is locally uniform in $v.$ $\square $

\bigskip

An alternative `representation' follows from Lemma 6.

\bigskip

\noindent \textbf{Theorem 6}$^{\prime }$\textbf{\ (Beurling Representation
Theorem).} \textit{For }$\varphi \in SN$\textit{\ with }$\varphi $ \textit{%
Baire/measurable eventually bounded away from 0, and }$f$\textit{\
measurable and }$\varphi $\textit{-regularly varying: there are }$\rho \in 
\mathbb{R}$\textit{, measurable }$d(.)\rightarrow d\in (0,\infty )$\textit{\
and continuous }$e(.)\rightarrow 0$ \textit{such that}%
\[
f(x)=d(x)\exp \left( \rho \int_{1}^{x}\frac{du}{\varphi (u)}%
+\int_{0}^{x}e(v)dv\right) =d(x)\exp \left( \rho \int_{1}^{x}\frac{u}{%
\varphi (u)}\frac{du}{u}+\int_{0}^{x}e(v)dv\right) ,
\]%
\textit{where for each }$t$%
\[
\int_{x}^{x+t\varphi (x)}e(v)dv=o(1).
\]

\noindent \textit{Proof.} Choose $\rho $ as in Theorem 3 and, referring to
the flow rate $\varphi (x)>0$ at $x,$ put%
\[
\tilde{h}(x):=h(x)-\rho \int_{1}^{x}\frac{du}{\varphi (u)}, 
\]%
where $h=\log f.$ Here, since $1/\varphi (x)$ is eventually bounded above as 
$x\rightarrow \infty $ and our analysis is asymptotic, w.l.o.g. we may
assume again by Luzin's Theorem that $\varphi $ here is continuous.

By Theorem 2 (UCT), locally uniformly in $t$ one has, as in Theorem 6, a
`reduction' formula for $\tilde{h}:$%
\begin{eqnarray}
&&\tilde{h}(x+t\varphi (x))-\tilde{h}(x)-\rho t+\rho \int_{x}^{x+t\varphi
(x)}\frac{du}{\varphi (u)} \\
&=&h(x+t\varphi (x))-h(x)-\rho t=o(1).  \label{re-1}
\end{eqnarray}

Fix $y$ and let $K>0$ be a bound for $1/\varphi $, far enough to the right.
We will use local uniformity in (\ref{re-1}) on the interval $|t|\leq
|y|K^{-1}.$ First, take $t=y/\varphi (x),$ so $|t|\leq |y|/K,$ so by (\ref%
{re-1}),%
\begin{eqnarray*}
\tilde{h}(x+y)-\tilde{h}(x) &=&\frac{\rho y}{\varphi (x)}-\rho \int_{x}^{x+y}%
\frac{du}{\varphi (u)}+o(1) \\
&=&\rho \int_{x}^{x+y}\left( \frac{1}{\varphi (x)}-\frac{1}{\varphi (u)}%
\right) du+o(1) \\
&=&\frac{\rho }{\varphi (x)}\int_{0}^{y}\left( 1-\frac{\varphi (x)}{\varphi
(x+w)}\right) d+o(1).
\end{eqnarray*}%
By Lemma 6 applied to the set $\{v:|v|\leq |y|\},$ which corresponds to the $%
w$ range in the integral above, we have 
\begin{equation}
\tilde{h}(x+y)-\tilde{h}(x)=\frac{\rho }{\varphi (x)}\int_{0}^{y}\left( 1-%
\frac{\varphi (x)}{\varphi (x+v)}\right) dv+o(1)=o(1),  \label{new-*}
\end{equation}%
since $1/\varphi (x)$ is bounded. That is, $\tilde{h}(x)$ is slowly varying
in the additive Karamata sense (as with log $\varphi $ in the lemma). So by
the Karamata-de Bruijn representation (see BGT, 1.3.3),%
\begin{equation}
\tilde{h}(x)=c(x)+\int_{0}^{x}e(v)dv,  \label{new-**}
\end{equation}%
for some measurable $c(.)\rightarrow c\in \mathbb{R}$ and continuous $%
e(.)\rightarrow 0.$ Re-arranging yields%
\[
\log f(x)=h(x)=\tilde{h}(x)+\rho \int_{1}^{x}\frac{du}{\varphi (u)}%
=c(x)+\rho \int_{1}^{x}\frac{du}{\varphi (u)}+\int_{0}^{x}e(v)dv.
\]

Taking $d(x)=e^{c(x)}$ we obtain the desired representation. To check this,
w.l.o.g. we now take $d(x)=1$, and continue by substituting $u=x+s\varphi
(x) $ to obtain from (\ref{new-*}) and (\ref{new-**}) 
\[
\int_{x}^{x+t\varphi (x)}e(v)dv=\tilde{h}(x+t\varphi (x))-\tilde{h}(x)=o(1), 
\]%
since $\varphi \in SN.$ $\square $

\bigskip

The proof above remains valid when $\rho =0$ for arbitrary $\varphi \in SN,$
irrespective of whether $\varphi $ is bounded away from zero or not. Since $%
\varphi $ is itself $\varphi $-regularly varying with corresponding index $%
\rho =0,$ we have an alternative to the Bloom-Shea representation of $%
\varphi $ via the de Bruijn-Karamata representation. We record this as

\bigskip

\noindent \textbf{Corollary 2.} \textit{For measurable }$\varphi \in SN$ 
\textit{there are measurable }$d(.)\rightarrow d\in (0,\infty )$\textit{\
and continuous }$e(.)\rightarrow 0$\textit{\ such that}%
\[
\varphi (x)=d(x)\exp \left( \int_{1}^{x}e(v)dv\right) ,
\]%
\textit{where for each }$t$%
\[
\int_{x}^{x+t\varphi (x)}e(v)dv=o(1).
\]

We have in the course of the proof of Theorem 4 in fact also shown:

\bigskip

\noindent \textbf{Corollary 3.} \textit{For }$\varphi \in SN$\textit{\
bounded below and }$f$\textit{\ measurable and }$\varphi $\textit{-regularly
varying: }%
\[
f(x)=\tilde{f}(x)\exp \left( \rho \int_{1}^{x}\frac{du}{\varphi (u)}\right)
, 
\]%
\textit{for some }$\rho \in \mathbb{R}$\textit{\ and some Karamata
(multiplicatively) slowly varying }$\tilde{f}.$

\bigskip

\noindent \textbf{7. Complements}\newline

\noindent 1. \textit{Direct and indirect specialization to the Karamata
framework. }That $\varphi (x)=o(x)$ is part of the definition of Beurling
slow variation (see \S 1). We note that allowing $\varphi (x)=x$ formally
puts us back in the framework of Karamata regular variation. The results
above all extend in this way, and in particular one formally recovers the
standard results of \textit{multiplicative} Karamata theory:%
\[
\varphi (x)=x,\qquad T_{t}^{\varphi }x=(t+1)x,\qquad f(x)=d(x)\exp \left(
\rho \int_{1}^{x}\frac{du}{u}+\int_{0}^{x}e(v)dv\right) .
\]%
Of course $\varphi (x)=1$ directly specializes Beurling to Karamata \textit{%
additive} regular variation; the latter is equivalent under log/exp
transformation to the Karamata multiplicative form, so in that sense the
Beurling theory above incorporates the whole of Karamata theory, albeit
indirectly.\newline
\noindent 2. \textit{Karamata- versus Beurling-variation in the case }$%
\varphi (x)=x$\textit{. }As the reader may verify, the proofs of Theorems
2-4 above apply also to the context of $\varphi (x)=x$ (amended around the
origin to satisfy $\varphi >0$ in Th. 4), where the restriction $\varphi
(x)=o(x)$ fails (that is, $\varphi \notin SN)$. However, for $\varphi (x)=x$
one has $\varphi (x+u\varphi (x))/\varphi (x)=1+u$ (locally uniformly), and
indeed this latter fact suffices for the Shift Lemma and for Theorem 2 (UCT)
as $y_{n}\rightarrow \infty $, and $c_{n}\rightarrow 1$ still hold. (This is
a matter of separate interest, addressed in full generality via the $\lambda 
$-UCT of [Ost2]; $\lambda (t)=1+t$ is the case relevant here.) Thus Theorems
2 and 3 above actually contain their Karamata multiplicative counterparts ($%
\varphi (x)=x)$ directly, rather than depending on them. This situation is
akin to that of the relationship between the de Haan theory of BGT, Ch. 3
and the Karamata theory of BGT, Ch. 1, cf. the `double-sweep' procedure of
BGT of p. 128, and \S\ 3.13.1 p. 188. Of course the exception is Theorem 6$%
^{\prime }$ where, though the Beurling case specializes to the Karamata
case, nevertheless it depends on the Karamata case, by design.\newline
\noindent 3. \textit{Miller homotopy. }Given Theorems 4 and 5, w.l.o.g. in
the asymptotic analysis one may take $\varphi \in \mathcal{C}^{1}$ --
continuously differentiable; then $H(x,t)=x+t\varphi (x)$ is a Miller
homotopy in the sense of [BinO9], i.e. its three defining properties are
satisfied: (i) $H(x,0)=x,$ (ii) $H_{x},H_{t}$ exist, and (iii) $%
H_{t}(t,x)=\varphi (x)>0$. In view of this, Miller's [Mil] generalization of
a result due to Kestelman and also Borwein and Ditor (see [BinO9] for their
context) asserts that \textit{for any non-negligible (Baire/measurable) set }%
$T$\textit{\ and any null sequence }$z_{n},$\textit{\ for quasi all }$t\in T$%
\textit{\ there is an infinite }$\mathbb{M}_{t}\subseteq \mathbb{N}$\textit{%
\ such that }$\{t+\varphi (t)z_{m}:m\in \mathbb{M}_{t}\}\subseteq T;$ cf.
Theorem 0 above.\newline
\noindent 4. \textit{Self-neglect in auxiliary functions. }We have seen in
the Beurling UCT\ how uniformity in the auxiliary function $\varphi $ passes
`out' to $\varphi $-regularly varying $f.$ For the converse (uniformity
passing `in' from $f$ to $\varphi $), we note the following.\newline
\newline
\noindent \textbf{Proposition.} \textit{For Baire (measurable) }$f,$\textit{%
\ if for some Baire (measurable) function }$\varphi >0$\textit{\ and some
real }$\rho \neq 0$%
\[
f(x+u\varphi (x))/f(x)\rightarrow e^{\rho u},\text{ locally uniformly in }u,
\]%
\textit{then }$\varphi \in SN.$

\bigskip

\noindent \textit{Proof. }Replace $u$ by $\rho u$ and $\varphi (x)$ by $\psi
(x)=\varphi (x)/\rho $ to yield 
\[
f(x+u\psi (x))/f(x)\rightarrow e^{u},
\]%
(allowing $u<0),$ then follow verbatim as in BGT, 3.10.6, which relies only
on uniformity (and on the condition that $x+u\psi (x)\rightarrow \infty $ as 
$x\rightarrow \infty ,$ which is also deduced from uniformity in BGT\
3.10.1). $\square $

\bigskip 

\begin{center}
\textbf{References}
\end{center}

\noindent \lbrack BajK] B. Baj\v{s}anski and J. Karamata, Regularly varying
functions and the principle of equi-continuity, \textsl{Publ. Ramanujan Inst.%
}, \textbf{1} (1968/69), 235-246.\newline
\noindent \lbrack Bal] A. A. Balkema, \textsl{Monotone transformations and
limit laws}. Mathematical Centre Tracts, No. 45. Mathematisch Centrum,
Amsterdam, 1973.\newline
\noindent \lbrack BalGdH] A. A. Balkema, J. L. Geluk and L. de Haan, An
extension of Karamata's Tauberian theorem and its connection with
complimentary convex functions. \textsl{Quart. J. Math. }\textbf{30} (1979),
385-416.\newline
\noindent \lbrack Ban] S. Banach, \textsl{Th\'{e}orie des operations lin\'{e}%
aire}. Reprinted in \textsl{Collected Works}, vol. II, 401-444, (PWN,
Warszawa 1979) (1st ed. 1932).\newline
\noindent \lbrack Bec] A. Beck, \textsl{Continuous flows on the plane},
Grundl. math. Wiss. \textbf{201}, Springer, 1974. \newline
\noindent \lbrack BinGT] N. H. Bingham, C. M. Goldie and J. L. Teugels, 
\textsl{Regular variation}, 2nd ed., Cambridge University Press, 1989 (1st
ed. 1987). \newline
\noindent \lbrack BinO1] N. H. Bingham and A. J. Ostaszewski, The index
theorem of topological regular variation and its applications. \textsl{J.
Math. Anal. Appl.} \textbf{358} (2009), 238-248. \newline
\noindent \lbrack BinO2] N. H. Bingham and A. J. Ostaszewski, Infinite
combinatorics and the foundations of regular variation. \textsl{J. Math.
Anal. Appl.} \textbf{360} (2009), 518-529. \newline
\noindent \lbrack BinO3] N. H. Bingham and A. J. Ostaszewski, Infinite
combinatorics in function spaces: category methods, \textsl{Publ. Inst.
Math. (Beograd)} (N.S.) \textbf{86 }(100) (2009), 55--73.\newline
\noindent \lbrack BinO4] N. H. Bingham and A. J. Ostaszewski, Beyond
Lebesgue and Baire II: bitopology and measure-category duality. \textsl{%
Colloq. Math.} \textbf{121} (2010), 225-238.\newline
\noindent \lbrack BinO5] N. H. Bingham and A. J. Ostaszewski, Topological
regular variation. I: Slow variation; II: The fundamental theorems; III:
Regular variation. \textsl{Topology and its Applications} \textbf{157}
(2010), 1999-2013, 2014-2023, 2024-2037. \newline
\noindent \lbrack BinO6] N. H. Bingham and A. J. Ostaszewski, Normed groups:
Dichotomy and duality. \textsl{Dissertationes Math.} \textbf{472} (2010),
138p. \newline
\noindent \lbrack BinO7] N. H. Bingham and A. J. Ostaszewski, Kingman,
category and combinatorics. \textsl{Probability and Mathematical Genetics}
(Sir John Kingman Festschrift, ed. N. H. Bingham and C. M. Goldie), 135-168,
London Math. Soc. Lecture Notes in Mathematics \textbf{378}, CUP, 2010. 
\newline
\noindent \lbrack BinO8] N. H. Bingham and A. J. Ostaszewski, Dichotomy and
infinite combinatorics: the theorems of Steinhaus and Ostrowski. \textsl{%
Math. Proc. Cambridge Phil. Soc.} \textbf{150} (2011), 1-22. \newline
\noindent \lbrack BinO9] N. H. Bingham and A. J. Ostaszewski, Homotopy and
the Kestelman-Borwein-Ditor Theorem, \textsl{Canadian Math. Bull.} 54.1
(2011), 12-20.\newline
\noindent \lbrack BinO10] N. H. Bingham and A. J. Ostaszewski, Steinhaus
theory and regular variation: De Bruijn and after. \textsl{Indagationes
Mathematicae} (N. G. de Bruijn Memorial Issue), to appear.\newline
\noindent \lbrack BinO11] N. H. Bingham and A. J. Ostaszewski, Uniformity
and self-neglecting functions, preprint (http://arxiv.org/abs/1301.5894).%
\newline
\noindent \lbrack BinO12] N. H. Bingham and A. J. Ostaszewski, Beurling
regular variation: III. Banach algebras, asymptotic actions and cocycles, in
preparation.\newline
\noindent \lbrack Blo] S. Bloom, A characterization of B-slowly varying
functions. \textsl{Proc. Amer. Math. Soc.} \textbf{54} (1976), 243-250. 
\newline
\noindent \lbrack Dom] J. Domsta: \textsl{Regularly varying solutions of
linear equations in a single variable: applications to regular iteration,}
Wyd. Uniw. Gda\'{n}skiego, 2002.\newline
\noindent \lbrack deH] L. de Haan, On regular variation and its applications
to the weak convergence of sample extremes, \textsl{Math. Centre Tracts} 
\textbf{32}, Amsterdam, 1970.\newline
\noindent \lbrack HewR] E. Hewitt and K. A. Ross, \textsl{Abstract harmonic
analysis, I Structure of topological groups, integration theory, group
representations}$.$ Grundl. math. Wiss. 115, Springer, 1963.\newline
\noindent \lbrack Kech] A. S. Kechris: \textsl{Classical Descriptive Set
Theory.} Grad. Texts in Math. 156, Springer, 1995.\newline
\noindent \lbrack Kor] J. Korevaar, \textsl{Tauberian theorems: A century of
development}. Grundl. math. Wiss. \textbf{329}, Springer, 2004. \newline
\noindent \lbrack Kucz] M. Kuczma, \textsl{An introduction to the theory of
functional equations and inequalities. Cauchy's equation and Jensen's
inequality.} 2nd ed., Birkh\"{a}user, 2009 [1st ed. PWN, Warszawa, 1985].%
\newline
\noindent \lbrack Kur] C. Kuratowski, \textsl{Topologie}, Monografie Mat. 20
(4th. ed.), PWN Warszawa 1958 [K. Kuratowski, \textsl{Topology}, Translated
by J. Jaworowski, Academic Press-PWN, 1966]. \newline
\noindent \lbrack Mas] A. Mas, Representation of gaussian small ball
probabilities in\textsl{\ }$\ell _{2},$ arXiv:0901.0264.\newline
\noindent \lbrack Meh] M. R. Mehdi, On convex functions, \textsl{J. London
Math. Soc.}, 39 (1964), 321-326.\newline
\noindent \lbrack Mil] H. I. Miller, \textsl{Generalization of a result of
Borwein and Ditor}, Proc. Amer. Math. Soc. 105 (1989), no. 4, 889--893.%
\newline
\noindent \lbrack Om] E. Omey, On the class gamma and related classes of
functions, preprint, 2011.
(http://www.edwardomey.com/pages/research/reports/2011.php).\newline
\noindent \lbrack Ost1] A. J. Ostaszewski, Regular variation, topological
dynamics, and the Uniform Boundedness Theorem, \textsl{Topology Proceedings}%
, \textbf{36} (2010), 305-336.\newline
\noindent \lbrack Ost2] A. J. Ostaszewski, Beurling regular equivariation
and the Beurling functional equation, preprint (2013)
(www.maths.lse.ac.uk/Personal/adam).\newline
\noindent \lbrack Oxt] J. C. Oxtoby, \textsl{Measure and category,} 2nd ed.,
Grad. Texts Math. \textbf{2}, Springer, 1980. \newline

\bigskip

\noindent Mathematics Department, Imperial College, London SW7 2AZ;
n.bingham@ic.ac.uk \newline
Mathematics Department, London School of Economics, Houghton Street, London
WC2A 2AE; A.J.Ostaszewski@lse.ac.uk\newpage 

\end{document}